\documentclass[review]{elsarticle}

\usepackage{graphicx,amssymb}
\usepackage{epsfig}

\usepackage{amsmath}    

\usepackage{framed}

\begin{document}

\def\eR{\mathbb{R}}

\begin{frontmatter}

\title{Geometrical properties of a class of systems with spiral trajectories in ${\eR}^3$ }

\author[fer]{L.\ Korkut}
\ead{luka.korkut@fer.hr}

\author[fer]{D.\ Vlah\corref{corresp}}
\ead{domagoj.vlah@fer.hr}

\author[fer]{V.\ \v Zupanovi\'c}
\ead{vesna.zupanovic@fer.hr}

\cortext[corresp]{Corresponding author. Tel. number +385 01 612 99 03}

\address[fer]{University of Zagreb, Faculty of Electrical Engineering and Computing, Unska 3, 10000 Zagreb, Croatia}

\journal{arXiv.org}

\begin{abstract}
Here we study a class of second-order nonautonomous differential equations, and the corresponding planar and spatial systems, from the point of view of fractal geometry. The fractal oscillatority of solutions at infinity is measured  by oscillatory and phase dimensions. The oscillatory dimension is defined as the box dimension of the reflected solution near the origin, while the phase dimension is defined as the box dimension of a trajectory of the corresponding planar system in the phase plane. Using the phase dimension of the second-order equation we compute the box dimension of a spiral trajectory of the spatial system, lying in  Lipschitzian or H\" olderian surfaces. This phase dimension of the second-order equation is connected to the asymptotics of the associated Poincar\'e map. Also, the box dimension of a trajectory of the reduced normal form with one eigenvalue equals to zero, and a pair of pure imaginary eigenvalues has been computed when limit cycles bifurcate from the origin.
\end{abstract}
\medskip




\end{frontmatter}

\newcommand{\kutija}[1]{\begin{framed}#1\end{framed}}

\newtheorem{theorem}{Theorem}
\newtheorem{cor}{Corollary}
\newtheorem{prop}{Proposition}
\newtheorem{lemma}{Lemma}
\newdefinition{defin}{Definition}
\newdefinition{remark}{Remark}

\font\csc=cmcsc10

\def\esssup{\mathop{\rm ess\,sup}}
\def\essinf{\mathop{\rm ess\,inf}}
\def\wo#1#2#3{W^{#1,#2}_0(#3)}
\def\w#1#2#3{W^{#1,#2}(#3)}
\def\wloc#1#2#3{W_{\scriptstyle loc}^{#1,#2}(#3)}
\def\osc{\mathop{\rm osc}}
\def\var{\mathop{\rm Var}}
\def\supp{\mathop{\rm supp}}
\def\Cap{{\rm Cap}}
\def\norma#1#2{\|#1\|_{#2}}

\def\C{\Gamma}

\let\text=\mbox

\catcode`\@=11
\let\ced=\c
\def\a{\alpha}
\def\b{\beta}
\def\c{\gamma}
\def\d{\delta}
\def\g{\lambda}
\def\o{\omega}
\def\q{\quad}
\def\n{\nabla}
\def\s{\sigma}
\def\div{\mathop{\rm div}}
\def\sing{{\rm Sing}\,}
\def\singg{{\rm Sing}_\ty\,}

\def\A{{\cal A}}
\def\F{{\cal F}}
\def\H{{\cal H}}
\def\W{{\bf W}}
\def\M{{\cal M}}
\def\N{{\cal N}}
\def\S{{\cal S}}

\def\ty{\infty}
\def\e{\varepsilon}
\def\f{\varphi}
\def\:{{\penalty10000\hbox{\kern1mm\rm:\kern1mm}\penalty10000}}
\def\ov#1{\overline{#1}}
\def\D{\Delta}
\def\O{\Omega}
\def\pa{\partial}

\def\st{\subset}
\def\stq{\subseteq}
\def\pd#1#2{\frac{\pa#1}{\pa#2}}
\def\sgn{{\rm sgn}\,}
\def\sp#1#2{\langle#1,#2\rangle}

\newcount\br@j
\br@j=0
\def\q{\quad}
\def\gg #1#2{\hat G_{#1}#2(x)}
\def\inty{\int_0^{\ty}}
\def\od#1#2{\frac{d#1}{d#2}}

\def\bg{\begin}
\def\eq{equation}
\def\bgeq{\bg{\eq}}
\def\endeq{\end{\eq}}
\def\bgeqnn{\bg{eqnarray*}}
\def\endeqnn{\end{eqnarray*}}
\def\bgeqn{\bg{eqnarray}}
\def\endeqn{\end{eqnarray}}

\def\bgeqq#1#2{\bgeqn\label{#1} #2\left\{\begin{array}{ll}}
\def\endeqq{\end{array}\right.\endeqn}

\def\abstract{\bgroup\leftskip=2\parindent\rightskip=2\parindent
        \noindent{\bf Abstract.\enspace}}
\def\endabstract{\par\egroup}

\def\udesno#1{\unskip\nobreak\hfil\penalty50\hskip1em\hbox{}
             \nobreak\hfil{#1\unskip\ignorespaces}
                 \parfillskip=\z@ \finalhyphendemerits=\z@\par
                 \parfillskip=0pt plus 1fil}
\catcode`\@=11

\def\cal{\mathcal}

\def\eN{\mathbb{N}}
\def\Ze{\mathbb{Z}}
\def\Qu{\mathbb{Q}}
\def\Ce{\mathbb{C}}

\def\osd{\mathrm{osd}\,}

\section{Introduction and motivation}

We found our mathematical inspiration in the book of  C. Tricot \cite{tricot}, where the author introduced a new approach for studying curves. He showed for some classes of smooth curves, nonrectifiable near the accumulation point, that fractal dimension called box dimension, can "measure" the density of accumulation. C. Tricot computed box dimension for class of spiral curves and chirps.
In this article by geometric properties of systems we mean type of solution curves, which are here spirals and chirps. Furthermore we distinguish rectifiable and nonrectifiable curves. Whereas box dimension of rectifiable curve is trivial, we proceed to study nonrectifiable curves using Tricot's fractal approach, and compute the box dimension.

Since 1970s dimension theory for dynamics has evolved into an independent field of mathematics. Together with Hausdorff dimension, box dimension was used to characterize dynamics, in particular chaotic dynamics having strange attractors, see \cite{fdd}.
We use the box dimension, because of countable stability of Hausdorff dimension, its value is trivial on all smooth nonrectifiable curves, while the box dimension is nontrivial, that is, larger than 1. The box dimension is  suitable tool for classification of nonrectifiable curves.
Analogously, box dimension is a good tool for analysis of discrete dynamical systems. Using box dimension we can study orbits of one-dimensional discrete system near its fixed point. Slow convergence to stable point means higher density of an orbit near its fixed point, which implies bigger box dimension. Fast convergence is related to trivial box dimension.

A natural idea is that higher density of orbits reveals higher multiplicity of the fixed point. The multiplicity of the fixed points is related to the bifurcations which could be produced by varying parameters of a given family of systems. Bifurcation theory provides a strategy for investigating the bifurcations that occur within a family.

In the paper \v Zubrini\'c and \v Zupanovi\'c \cite{zuzu}, the number of limit cycles which could be produced from  weak foci and limit cycles is directly related to the box dimension of any trajectory. It was discovered that the box dimension of a spiral trajectory of weak focus signals a moment of Hopf and Hopf-Takens bifurcation. The result was obtained using Takens normal form \cite{takens}. Using numerical algorithm for computation of box dimension of trajectory, it is possible to predict change of stability of the system, through Hopf bifurcation.  Recent results Marde\v si\' c,  Resman and \v Zupanovi\' c  \cite{MRZ}, Resman \cite{majaformal}, Horvat Dmitrovi\' c  \cite{laho}, and Resman \cite{majaabel}, show efficiency of this approach to the bifurcation theory. From asymptotic expansion of $\varepsilon$-neighborhood of an orbit, we  read box dimension and Minkowski content from the leading term. In the mentioned articles, it has been showed that more information about dynamical system
could be read from other terms of the asymptotic expansion of $\varepsilon$-neighborhood.

In this article we study nonautonomous differential equation of second order, and the corresponding systems  with spiral trajectories, in $\eR^2$ and $\eR^3$. The planar system has the same type of spiral as in Takens normal form, see \cite{takens}, which is spiral with analytic first return Poincar\' e map, also having the same asymptotics in each direction. Here we studied graph of solution of differential equation, as well as the corresponding trajectories in the phase plane.  The system could produce limit cycles under perturbation, but it is left for further research.

According to the idea of qualitative theory of differential equations, oscillations of a class of second-order differential equations have been considered by phase plane analysis, in Pa\v si\'c, \v Zubrini\'c and \v Zupanovi\'c  \cite{chirp}. The novelty was a fractal approach, connecting the box dimension of the graph of a solution and the box dimension of a trajectory in the phase plane. Oscillatory and phase dimensions for a class of second-order differential equations have been introduced.
The notion of fractal oscillatority near a point for real functions of a real variable has been introduced in Pa\v si\'c \cite{pasic}, and studied for  second-order linear differential equations of Euler type, as the basic model. The study was based on the fact that the Euler  equation has chirp-like oscillatory solutions.
This fractal approach to the standard Bessel equation and a generalized Bessel equation, can be found in Pa\v si\'c and  Tanaka \cite{pasic_bessel_jmaa,mersat},  and  also in \cite{bessel}.

On the other hand, in  \cite{zuzu,zuzulien}, the box dimension of spiral trajectories of a system with pure imaginary eigenvalues, near singular points and limit cycles, has been studied using normal forms, and the Poincar\' e map. These results are based on the fact that these spiral trajectories are  of power type in polar coordinates.
Also, it is shown that the box dimension is sensitive with respect to bifurcations, e.g.,  it jumps from the trivial value $1$ to the value $4/3$  when the Hopf bifurcation occurs. Degenerate Hopf bifurcation or Hopf-Takens bifurcation can reach even larger box dimension of a  trajectory near a singular point. This value is related to the multiplicity of the singular point surrounded with spiral trajectories. This phenomenon has been discovered for  discrete systems in  Elezovi\'c, \v Zupanovi\'c and \v Zubrini\'c \cite{neveda} concerning saddle-node and period-doubling bifurcations, and generalized in  \cite{laho}.
Also, in  \cite{MRZ}   there are  results about multiplicity of the  Poincar\' e map   near a weak focus, limit cycle, and saddle-loop, obtained using the box dimension. Isochronicity of a focus has been characterized by box dimension in Li and Wu \cite{li}. Formal normal forms for parabolic diffeomorphisms have been characterized by fractal invariants of the $\varepsilon-$neighborhood of a discrete orbit in \cite{majaformal}.  All these results are related to the 16th Hilbert problem.

This work is a part of our research, which has been undertaken in order to understand  relation between the graph of certain type of oscillatory  function, and the corresponding spiral curve in the phase plane. We believed that chirp-like oscillations  defined by $X(\tau)=\tau^{\a}\sin1/\tau$ ``correspond" to spiral oscillations  $r=\varphi ^{-\a}$ in the phase plane, in polar coordinates.
The relation  between these two objects has been established in \cite{cswavy}, introducing a new notion of the wavy spiral. Applications include two directions. Roughly speaking, we consider spirals generated by chirps, and chirps generated by spirals. If we know behavior in the phase plane, we can obtain the behavior of the corresponding graph, and vice versa.
As examples we may consider weak foci of planar autonomous systems, including the Li\' enard equation, because all these singularities are of spiral power type $r=\varphi ^{-\a}$, $\a\in(0,1)$, see \cite{chirp,zuzu}. As an application of the converse direction, from a chirp to the spiral,  we were looking for the second-order equation with an oscillatory solution having chirp-like behavior. The Bessel equation of order $\nu$ is a nice example of similar behavior, see \cite{bessel}. Whereas the Bessel equation is a second-order nonautonomous equation, we interpret the equation as a system in $\eR^3$, using $t\to\infty$, instead of the standard approach with a variable near the origin. The system studied in this article, see (\ref{alfap}), coincides with the Bessel system for $p(t)=t^{-\alpha}$, $\alpha=\nu=1/2$, and $q(t)=t$. We classify trajectories of the system with respect to their geometrical and fractal properties.

Why we study functions which behave like $X(\tau)=\tau^{\a}\sin1/\tau$, and $r=\varphi ^{-\a}$ in polar coordinates?
Our starting point is Tricot's book which gives us formulas for box dimension of $X(\tau)=\tau^{\a}\sin1/\tau^\beta$, $0<\alpha<\beta<1$, and $r=\varphi ^{-\a}$ , $0<\alpha<1$. For other parameters $\alpha$, $\beta$  these curves are rectifiable. We wanted to analyze power spirals which have same asymptotics of the  Poincar\'e map in all directions. Poincar\'e map which corresponds to weak focus is analytic, and limit cycles bifurcate in the classical Hopf bifurcation. Poincar\'e  maps near general foci, nilpotent or degenerate, as well as near polycycles are not analytic and the logarithmic terms show up in the asymptotic expansion, see Medvedeva \cite{medvedeva} and Roussarie \cite{16th-hilbert}. In that case Poincar\'e  map has different asymptotics, showing characteristic directions,  see Han, Romanovski  \cite{hanrom}. Nilpotent focus has two different asymptotics, so we can relate that focus with two chirps with different asymptotics.
Here, we study foci related to one chirp. Why we have chirps with $\beta=1$? For $\alpha+1\le\beta$ we have curves which do not accumulate in the origin, while for $\alpha+1>\beta$, if $\beta\ne1$ it is easy to see that the spiral converges to zero in ``oscillating'' way. Wavy spirals appear in that situation, see \cite{bessel}, \cite{cswavy}. Curves which are spirals with self intersections like springs, could be defined by oscillatory integrals, so they appear as a generalization of the clothoid defined by Fresnel integrals. Asymptotics of the oscillatory integrals, which are related to singularity theory, could be found in Arnold \cite{arnold-vol2}. Their fractal analysis is our work in progress.
Furthermore, fixing $\b=1$ we achieve the whole interval of nonrectifiability both for spirals and corresponding chirps.

Also, the results about spiral trajectories in $\eR^3$, from \v Zubrini\'c and \v Zu\-pa\-no\-vi\'c  \cite{cras,bilip} are extended to the systems where some kind of Hopf bifurcation occurs.  The box dimension of a trajectory of the reduced normal form with one zero eigenvalue, and a pair of pure imaginary eigenvalues, has been computed at the moment of the birth of a limit cycle. Essentially, the Hopf bifurcation studied here is a planar bifurcation, but the third equation affects the box dimension of the corresponding trajectory in the space. We show that in $3$-dimensional space, a limit cycle bifurcates with the box dimension of a spiral trajectory larger than $4/3$, which is the value of the standard planar Hopf bifurcation.

\smallskip

Our intention is to understand a fractal connection between oscillatority of solutions of differential equations and oscillatority of their trajectories in the phase space. Our work is mostly motivated by two nice formulas from the monograph of C. Tricot \cite[p.\ 121]{tricot}. He computed the box dimension for a class of chirps and for a class of spirals of power type in polar coordinates. We are looking for a model to apply these formulas, and also to show that chirps and spirals are a different manifestations of the same phenomenon.
Here we study, as a model a class of second-order nonautonomous equations, exhibiting both chirp and spiral behavior

\bgeq\label{eqalfap}
\ddot x - \left[\frac{2\,p'(t)}{p(t)} + \frac{q''(t)}{q'(t)}\right]\dot x + \left[q'^2(t) + \frac{2\,p'^2(t)}{p^2(t)} - \frac{p''(t)}{p(t)} + \frac{p'(t)q''(t)}{p(t)q'(t)}\right]x = 0,
\endeq
$t\in [t_0,\infty)$, $t_0>0$, where $p$ and $q$ are functions of class $C^2$.
The explicit solution is $x(t)=C_1p(t)\sin q(t) +C_2p(t)\cos q(t) $, which is a chirp-like function. If $z=(\gamma/(t-C_3))^{\gamma}$, $\gamma>0$, we get the cubic system

\bgeqn\label{alfap}
\dot x&=&y \nonumber  \\
\dot y&=& -U(z) x +V(z) y  \\
\dot z&=& -z^{\delta},\quad  z\in (0,z_0] , \nonumber
\endeqn
where $\delta:=(\gamma+1)/\gamma > 1$ and
\bgeqn
U(z)&:=&q'^2(\gamma z^{-\frac{1}{\gamma}}) + \frac{2\,p'^2(\gamma z^{-\frac{1}{\gamma}})}{p^2(\gamma z^{-\frac{1}{\gamma}})} - \frac{p''(\gamma z^{-\frac{1}{\gamma}})}{p(\gamma z^{-\frac{1}{\gamma}})} + \frac{p'(\gamma z^{-\frac{1}{\gamma}})q''(\gamma z^{-\frac{1}{\gamma}})}{p(\gamma z^{-\frac{1}{\gamma}})q'(\gamma z^{-\frac{1}{\gamma}})} ,\nonumber\\
V(z)&:=&\frac{2\,p'(\gamma z^{-\frac{1}{\gamma}})}{p(\gamma z^{-\frac{1}{\gamma}})} + \frac{q''(\gamma z^{-\frac{1}{\gamma}})}{q'(\gamma z^{-\frac{1}{\gamma}})} .\nonumber
\endeqn
It has a spiral trajectory in $\mathbb{R}^3$. In the special case $\gamma=1$, we get $z=1/(t-C_3)$ and $\delta=2$.

In this article we compute the box dimension of a spiral trajectory of the system (\ref{alfap}) exploiting the dimension of $(\a,1)-$chirp $X(\tau)=\tau^{\alpha}\sin 1/\tau$, $\a\in(0,1)$, for $\tau>0$ small, and also the dimension of the wavy spiral, see \cite{cswavy}. Using a change of variable for time variable $\tau\mapsto\tau^{-1}$, the  infinity is mapped to the origin, and such reflected  solution of (\ref{eqalfap}) with respect to time is called  the {\it reflected solution}.
We use notation $t=\tau^{-1}$. If function $p(t)$ in (\ref{eqalfap})  is ``similar'' to $t^{-\a}$, and function $q(t)$ is ``similar'' to $t$, then the reflected solution of $x(t)$
\bgeqn
X(\tau)&=&C_1p(\frac{1}{\tau})\sin q(\tau^{-1})+ C_2p(\frac{1}{\tau})\cos q(\tau^{-1})\nonumber \\
&=&\sqrt{C_1^2+C_2^2}\ p(\frac{1}{\tau})\sin(q(\tau^{-1})+\arctan{\frac{C_2}{C_1}}), \quad \tau\in(0,\frac{1}{t_0}],\nonumber
\endeqn
is an $(\a,1)$-chirp-like function near the origin, see \cite{lukamaja}.
Before we obtained results connecting functions ``similar'' to $(\a,1)-$chirps, and  spirals ``similar'' to $r=\varphi ^{-\a}$, $\a\in(0,1)$, in the phase plane, see \cite{cswavy}. Applications include nonautonomous planar systems, so here we introduce the third variable $z$ depending on the time $t$. Furthermore, the box dimension of a  trajectory depends on $\gamma>0$. For some values of $\gamma$ trajectory in $\eR ^3$ is obtained as bi-Lipschitzian image of the spiral from the phase plane, which does not affect the box dimension. For other values, trajectory lies in the H\" olderian surface, affecting the box dimension. The H\" olderian surface has an infinite derivative at the origin, which is the point of accumulation of the spiral. Spirals of the H\" olderian type have the ``tornado shape'' with  a small bottom and wide top.

It is interesting to notice that our results  about the box dimension of  planar trajectories of a system with pure imaginary eigenvalues, show that the box dimension of any trajectory depends on the exponents of the system. In $\eR ^3$ we have already found an example,  see \cite{cras,bilip}, where dimension depends on the coefficients of the systems, which will be the case in (\ref{alfap}).
See \cite{bilip} for the computation of the box dimension of  the system
\bgeq\label{spp}
\dot r=a_1rz,\q
\dot \f=1,\q
\dot z=b_2z^2,
\endeq
in cylindrical coordinates. If $a_1/b_2\in(0,1]$ then any spiral trajectory $\Gamma$ of (\ref{spp}) has the box dimension $\dim_B\C=\frac{2}{1+a_1/b_2}$ near the origin.

\smallskip

\section{Definitions}

Let us introduce some definitions and notation.
For $A\st\eR^N$  bounded
we define
 {\it{  $\e$-neighborhood}} of  $A$ as:
$
A_\e:=\{y\in\eR^N\:d(y,A)<\e\}
$.
By {\it lower $s$-dimensional  Minkowski content} of $A$, $s\ge0$ we mean
$$
\M_*^s(A):=\liminf_{\e\to0}\frac{|A_\e|}{\e^{N-s}},
$$
and analogously for the {\it upper $s$-dimensional Minkowski content} $\M^{*s}(A)$.
The lower and upper box dimensions of $A$ are
$$
\underline\dim_BA:=\inf\{s\ge0\:\M_*^s(A)=0\}
$$
 and analogously
$\ov\dim_BA:=\inf\{s\ge0\:\M^{*s}(A)=0\}$.
If these two values coincide, we call it simply the box dimension of $A$, and denote by $\dim_BA$. It will be our situation.
If $0<\M_*^d(A)\le\M^{*d}(A)<\ty$ for some $d$, then we say
 that $A$ is {\it Minkowski nondegenerate}. In this case obviously $d=\dim_BA$.
In the case when lower or upper $d$-dimensional Minkowski contents of $A$ are $0$ or $\ty$,
where $d=\dim_BA$,
we say that $A$ is {\it degenerate}.
For more details on these definitions see e.g.\ Falconer \cite{falc}, and \cite{zuzu}.
\smallskip

Let $x:[t_0,\ty)\to\eR$, $t_0>0$, be a continuous function. We say that $x$ is an {\it oscillatory function}
 near $t=\ty$
if there exists a sequence $t_k\searrow \ty$
such that $x(t_k)=0$, and functions $x|_{(t_k,t_{k+1})}$ intermittently change sign for $k\in \eN$.

\smallskip

Let  $u:(0,t_0]\to\eR$, $t_0>0$, be a continuous function. We say that $u$ is an {\it oscillatory function} near the origin if there exists a sequence $s_k$ such that $s_k\searrow 0$ as $k\rightarrow \infty $, $u(s_k)=0$ and restrictions $u|_{(s_{k+1},s_k)}$ intermittently change sign, $k\in \eN$.

\smallskip

Let us define $X:(0,1/t_0]\to\eR$ by $X(\tau)=x(1/\tau)$. We say that $X(\tau)$ is oscillatory near the origin if $x=x(t)$ is oscillatory near $t=\ty$.
We measure the rate of oscillatority of $x(t)$ near $t=\ty$ by the rate of oscillatority of $X(\tau)$ near $\tau=0$.
More precisely, the {\it oscillatory dimension} $\dim_{osc}(x)$ (near $t=\ty$) is defined as the box dimension
of the graph of $X(\tau)$ near $\tau=0$. In Radunovi\' c, \v Zubrini\'c and \v Zupanovi\'c \cite{RaZuZu} box dimension of unbounded sets has been studied.

\smallskip

Assume now that $x$ is of class $C^1$. We say that $x$ is a {\it phase oscillatory} function if the following stronger condition holds: the set
$\C=\{(x(t),\dot x(t)):t\in[t_0,\ty)\}$ in the plane is a spiral converging to the origin.

By a spiral here we mean the graph of a function $r=f(\varphi)$, $\varphi\geq\varphi_1>0$, in polar coordinates, where
\begin{equation}\label{def_spiral}
\left\{\begin{array}{l}
f:[\varphi_1,\infty)\rightarrow(0,\infty) \textrm{ is such that } f(\varphi)\to 0 \textrm{ as } \varphi\to\infty,\\
f \textrm{ is \emph{radially decreasing} (i.e., for any fixed } \varphi\geq\varphi_1\\
\textrm{the function } \mathbb{N}\ni k\mapsto f(\varphi+2k\pi) \textrm{ is decreasing)},
\end{array}\right. \nonumber
\end{equation}
which is the  definition from \cite{zuzu}. By a spiral we also mean a mirror image of the spiral (\ref{def_spiral}), with respect to the $x$-axis.

The {\it phase dimension} $\dim_{ph}(x)$ of the function $x(t)$ is defined as the box dimension of the corresponding spiral $\C=\{(x(t),\dot x(t)):t\in[t_0,\ty)\}$.

\smallskip
We use a result for box dimension of  graph $G(X)$ of standard {\it $(\a,\b)$-chirps}
defined by
\bgeq\label{standard_chirp}\nonumber
X_{\a,\b}(\tau)=\tau^{\a}\sin(\tau^{-\beta})\nonumber.
\endeq
For $0<\a<\b$ we have
\bgeq\label{dim}\nonumber
\dim_BG(X_{\a,\b})=2-(\a+1)/(\b+1)\nonumber,
\endeq
 and the same for $X_{\a,\b}(\tau)=\tau^{\a}\cos(\tau^{-\beta})$,
see Tricot  \cite[p.\ 121]{tricot}.
Also we use a result for box dimension of spiral $\C$
defined by $r=\f^{-\a}$, $\f\ge\f_0>0$,  $\dim_B\C=2/(1+\a)$ when $0<\a\le1$, see Tricot \cite[p.\ 121]{tricot} and some generalizations from \cite{zuzu}.

Oscillatory and phase dimensions are fractal dimensions, which are well known
tool in the study of dynamics, see survey article \cite{fdd}.

\smallskip

For two real functions $f(t)$ and $g(t)$ of real variable
we write $f(t)\simeq g(t)$, and say that functions are {\it comparable} as $t\to0$ (as $t\to\ty$),
if there exist positive constants $C$ and $D$ such that $C\,f(t)\le g(t)\le D\,f(t)$ for all $t$ sufficiently close to $t=0$ (for all $t$ sufficiently large). For example, for a function $F:U\to V$ with $U,V\st\eR^2$, $V=F(U)$, the condition $|F(t_1)-F(t_2)|\simeq
|t_1-t_2|$ means that $f$ is a bi-Lipschitz mapping, i.e.,\ both $F$ and $F^{-1}$ are Lipschitzian.

We say that {\it function $f$ is comparable of class $k$ to  power $t^{-\alpha}$} if $f$ is class $C^k$ function, and $f^{(j)}(t)\simeq t^{-\a-j}$ as $t\to\ty$, $j=0,1,2,\dots, k.$

\smallskip

Also, we write $f(t)\sim g(t)$ if $f(t)/g(t)\to1$ as $t\to\infty$, and say that {\it function $f$ is comparable of class $k$ to  power $t^{-\alpha}$ in the limit sense}  if $f$ is class $C^k$ function, $f(t)\sim t^{-\a}$ and  $f^{(j)}(t)\sim (-1)^j\alpha(\alpha+1)(\alpha+j-1) t^{-\a-j}$ as $t\to\ty$, $j=1,2,\dots, k.$

We write $f(t)=O(g(t))$ as $t\to0$ (as $t\to\ty$) if there exists positive constant $C$ such that $|f(t)|\leq C|g(t)|$. We write $f(t)=o(g(t))$ as $t\to\ty$ if for every positive constant $\varepsilon$ it holds $|f(t)|\leq \varepsilon|g(t)|$ for all $t$ sufficiently large.

\smallskip

In the sequel we shall consider the functions of the form $y=p(\tau)\sin(q(\tau))$ or $y=p(\tau)\cos(q(\tau))$. If $p(\tau)\simeq \tau^\a$, $q(\tau)\simeq \tau^{-\b}$, $q'(\tau)\simeq \tau^{-\b-1}$ as $\tau\to 0$ then we say that $y$ is an $(\a,\b)$-{\it chirp-like function}.

\smallskip

\section{Spiral trajectories in $\eR^3$}

In this section we describe solutions of equation (\ref{eqalfap}) and trajectories of system (\ref{alfap}), with respect to box dimension, specifying a class of functions $p$ and $q$. Let $p(t)$ be comparable to power $t^{-\alpha}$, $\alpha>0$, in the limit sense, and let $q(t)$ be comparable to $K t$, $K>0$, in the limit sense.
Depending on ${\alpha}$, we have rectifiable spirals with trivial box dimension equal to $1$, or nonrectifiable spirals with nontrivial box dimension greater than $1$.
The box dimension will not exceed $2$ even in $\eR^3$, because these spirals lie on a surface. Mapping spiral from the plane to the Lipschitzian surface does not affect the box dimension, see \cite{cras,bilip}, while mapping to H\"olderian surface affects the box dimension.

In order to explain fractal behavior of the system (\ref{alfap}) we need a lemma  dealing with a bi-Lipschitz map. The idea is to use a generalization of the result about box dimension of a class of planar spirals from \cite[Theorem~4]{cswavy}, see Theorem \ref{gen_wavy} in the Appendix. It is well known result from \cite{falc}, that box dimension is preserved by bi-Lipschitz map. Putting together these two results we will obtain desired results about (\ref{alfap}).
For the sake of simplicity, we deal with trajectory $\Gamma$ of the solution of the system $(\ref{alfap})$  defined by
\bgeqn\label{spiral}
 x(t)&=&p(t) \sin q(t) \nonumber  \\
 y(t)&=& p'(t)\sin q(t) +p(t)q'(t)\cos q(t)   \\
 z(t)&=&\frac{1}{t^{\gamma}} .\nonumber
\endeqn
We can assume, without the loss of generality, that $q(t)$ is comparable to $t$, in the limit sense, by contracting time variable $t$ by factor $K$ and also contracting $x$ by factor $K^{\alpha}$, $y$ by factor $K^{\alpha+1}$, and $z$ by factor $K^{\gamma}$. Notice that rescaling of spatial variables by a constant factor is a bi-Lipschitz map, so the box dimension of trajectory $\Gamma$ is preserved.

Trajectory  $\Gamma$  has projection $\Gamma_{xy}$ to $(x,y)$-plane which is a planar spiral satisfying conditions of Theorem \ref{gen_wavy}, see the Appendix. In the following lemma we will prove that the mapping between planar spiral $\Gamma_{xy}$ and spacial spiral $\Gamma$ is bi-Lipschitzian near the origin. We prove lemma using definition of bi-Lipschitz mapping. Interesting phenomenon appeared in spiral $\Gamma_{xy}$, defined in polar coordinates and generated by a chirp. The radius $r(\varphi)$  is not decreasing function, there are some regions where $r(\varphi)$  increases causing some waves on the spiral. We introduced notion of wavy spiral in \cite{cswavy}. Also, the waves are found in the spiral generated by Bessel functions, and  by generalized Bessel functions, depending on the parameters in the equation, see \cite{bessel}. Furthermore, the surface containing the space spiral $\Gamma$ contains points with infinite derivative, showing  some vertical  regions.

\begin{lemma}\label{Domagoj}
Let the map $B:\eR^2\times\{0\}\to \eR^3$ be defined as $B(x(t),y(t),0)=(x(t),y(t),z(t))$, where coordinate functions are given by (\ref{spiral}). Let $p(t)\in C^2$ is comparable of class $1$  to $t^{-\a}$, $\a\in(0,1)$, in the limit sense, and $p''(t)\in o(t^{-\alpha})$, as $t\to\infty$. Let $q(t)\in C^2$ is comparable of class $1$  to $K t$, $K>0$, in the limit sense, and $q''(t)\in o(t^{-2})$, as $t\to\infty$.
Let $\Gamma$ is defined by parametrization $(x(t),y(t),z(t))$ from (\ref{spiral}) and $\Gamma_{xy}$ is the projection of $\Gamma$ to $(x,y)$-plane.
If $\gamma\geq\alpha$ then map $B|_{\Gamma_{xy}}$ is bi-Lipschitzian near the origin.
\end{lemma}
{\it Proof.}
Without the loss of generality we assume $K=1$.

It is clear that $B(\Gamma_{xy})=\Gamma$. We have to prove that there exist two positive constants $K_1,K_2$ such that
\bgeqn\label{bil}
K_1 d( (x(t_1),y(t_1),0), (x(t_2),y(t_2),0))\le \\
d( B(x(t_1),y(t_1),0), B(x(t_2),y(t_2),0))\le  \nonumber\\
K_2 d( (x(t_1),y(t_1),0), (x(t_2),y(t_2),0)) \nonumber,
\endeqn
where $d$ is Euclidian metrics and $t_1,t_2>t_0$, for $t_0$ sufficiently large. Notice that, without loss of generality, $t_1\leq t_2$. It is obvious that by $K_1=1$ the left hand side inequality is satisfied.
In order to prove right hand side inequality, first  we prove
\bgeq
(z(t_1)-z(t_2))^2\le C\left((x(t_1)-x(t_2))^2+(y(t_1)-y(t_2))^2\right), \nonumber
\endeq
and then the right inequality will be satisfied. From the proof of the planar case, Theorem \ref{gen_wavy}, we know that
\bgeq\label{fi}
\f (t) =t+\frac{\pi}2 + O(t^{-1})=\frac{\pi}2+t\left(1+O\left(t^{-2}\right)\right), \quad t\to \infty \nonumber
\endeq
\bgeq\label{fialfa}
f(\f)\simeq \f ^{-\alpha}, \quad \varphi\to \infty .
\endeq
From generalization of \cite[Lemma 3]{cswavy} used in the proof of Theorem \ref{gen_wavy}, using assumptions on $p$ and $q$, it follows that there exists $C_1\in (0,1)$, such that for every $\Delta\f$, $\frac{\pi}{3}\le \Delta\f\le 2\pi+\frac{\pi}{3}$, holds
$$
f(\f)-f(\f + \Delta\f )\ge  \Delta\f \alpha C_1\f ^{-\alpha-1},
$$
for $\f$ sufficiently large. Let
\bgeqn\label{fi12}
\f_1=\f(t_1)=\pi/2 +t_1\left(1+O\left(t_1^{-2}\right)\right),\\
\f_2=\f(t_2)=\pi/2 +t_2\left(1+O\left(t_2^{-2}\right)\right). \nonumber
\endeqn

We first consider several cases where $\alpha\leq\gamma<1$. First, let $|\f_2-\f_1|\le \frac{\pi}3$.
From \cite[Proposition 1]{bessel}, and (\ref{fialfa}), (\ref{fi12}), we have
\bgeqn
\sqrt{(x(t_1)-x(t_2))^2+(y(t_1)-y(t_2))^2}\ge \frac{2}{\pi}(\f_2-\f_1)\min\{f(\f_1),f(\f_2)\}\ge \nonumber\\
C_2(t_2-t_1)t_2^{-\alpha}. \nonumber
\endeqn
Hence,
\bgeqn\label{b}
(z(t_1)-z(t_2))^2=\left(\frac1{t_1^{\gamma}}-\frac1{t_2^{\gamma}}\right)^2=\frac{(t_2^{\gamma}-t_1^{\gamma})^2}{t_1^{2\gamma}t_2^{2\gamma}}\leq \frac{(t_2-t_1)^2(t_2^{\gamma-1}+t_1^{\gamma-1})^2}{t_1^{2\gamma}t_2^{2\gamma}} \leq \\
c_2\frac{(t_2-t_1)^2 t_1^{2(\gamma-1)}}{t_1^{2\gamma}t_2^{2\gamma}} = c_2\frac{(t_2-t_1)^2}{t_1^2t_2^{2\gamma}}\cdot\frac{C_2^2t_2^{-2\alpha}}{C_2^2t_2^{-2\alpha}}\le \nonumber\\
c_2\frac{(x(t_1)-x(t_2))^2+(y(t_1)-y(t_2))^2}{ C_2^2 t_1^2t_2^{2(\gamma-\alpha)}}\le \nonumber\\
C\left((x(t_1)-x(t_2))^2+(y(t_1)-y(t_2))^2\right). \nonumber
\endeqn
For the case $2\pi + \frac{\pi}3\ge |\f_2-\f_1|\ge \frac{\pi}3$, we have
\bgeqn
\sqrt{ (x(t_1)-x(t_2))^2+(y(t_1)-y(t_2))^2 }\ge f(\f_1)-f(\f_2)= \nonumber\\
 f(\f_1)-f(\f_1+(\f_2-\f_1))\ge C_1 ( \f_2-\f_1)\alpha{\f_1}^{-\alpha -1}\ge
 C_3 t_1^{-\alpha -1}(t_2-t_1). \nonumber
\endeqn
Then again
\bgeqn\label{c}
(z(t_1)-z(t_2))^2=c_3\frac{(t_2-t_1)^2}{t_1^2t_2^{2\gamma}}\le
c_3\frac{ (x(t_1)-x(t_2))^2+(y(t_1)-y(t_2))^2 }{C_3^2 t_1^{2(\gamma-\alpha)}}\le \\
C \left( (x(t_1)-x(t_2))^2+((y(t_1)-y(t_2))^2 \right)\nonumber.
\endeqn
For the case $|\f_2-\f_1|\ge 2\pi + \frac{\pi}3$, we define $n:=\left[\frac{\f_2-\f_1- \frac{\pi}3}{2\pi}\right]$. Then we have
\bgeqn
\sqrt{(x(t_1)-x(t_2))^2+(y(t_1)-y(t_2))^2 }\ge f(\f_1)-f(\f_2)= \nonumber\\
\sum_{i=0}^{n-1}\left(f(\f_1+2i\pi)- f(\f_1+(i+1)2\pi)\right) + f(\f_1+2n\pi)-f(\f_2)\ge \nonumber\\
\sum_{i=0}^{n-1}2\pi\alpha C_1(\f_1+2i\pi)^{-\alpha-1}+\frac{\pi}3\alpha C_1(\f_1+2n\pi)^{-\alpha-1}\ge \nonumber\\
\frac{\pi}3\alpha C_1\sum_{i=0}^{n}(\f_1+2i\pi)^{-\alpha-1}= \frac{\pi}3\alpha C_1\sum_{i=0}^{n}{(2\pi)}^{ -\alpha-1}(\frac{\f_1}{2\pi}+i)^{ -\alpha-1} \ge \nonumber \\
\frac{\pi}3\alpha C_1{2\pi}^{ -\alpha-1}\int_{\frac{\f_1}{2\pi}}^{\frac{\f_1}{2\pi}+n-1}x^{ -\alpha-1}dx\ge C_4\f_1^{-\alpha-1}(\f_2-\f_1)\ge \nonumber \\ C_5 t_1^{-\alpha-1}(t_2-t_1).\nonumber
\endeqn
Furthermore
\bgeqn\label{d}
(z(t_1)-z(t_2))^2=c_5\frac{(t_2-t_1)^2}{t_1^2t_2^{2\gamma}}\le
c_5\frac{ (x(t_1)-x(t_2))^2+(y(t_1)-y(t_2))^2 }{C_5^2 t_1^{2(\gamma-\alpha)}}\le \\
C\left((x(t_1)-x(t_2))^2+(y(t_1)-y(t_2))^2\right) \nonumber.
\endeqn
From (\ref{b}), (\ref{c}), (\ref{d}), the right hand side inequality (\ref{bil}) follows with $K_2=\sqrt{1+C}$, where $C$ is (in all three cases) sufficiently small if $t_0$ is large enough.

For $t_0$ sufficiently large, it is easy to see that
$$
(z(t_1)-z(t_2))^2=\left(\frac{1}{t_1^{\gamma}}-\frac{1}{t_2^{\gamma}}\right)^2 \leq \left(\frac{1}{t_1^2}-\frac{1}{t_2^2}\right)^2 ,
$$
if $\gamma>2$. On the other hand, for $1\leq\gamma\leq 2$, considering
$$
(z(t_1)-z(t_2))^2\leq \frac{(t_2-t_1)^2(t_2^{\gamma-1}+t_1^{\gamma-1})^2}{t_1^{2\gamma}t_2^{2\gamma}} \leq
c_6\frac{(t_2-t_1)^2 t_2^{2(\gamma-1)}}{t_1^{2\gamma}t_2^{2\gamma}} = c_6\frac{(t_2-t_1)^2}{t_1^{2\gamma}t_2^2} ,
$$
the rest of the proof is analogous as for the case $\alpha\leq\gamma<1$.
\qed

\begin{theorem}{\label{phase}} {\rm (Trajectory in $\eR^3$)} Let $p(t)\in C^3$ be a function
 comparable of class $2$ to  power $t^{-\alpha}$, $\alpha>0$, in the limit sense, and $p^{(3)}(t)\in O(t^{-\alpha-3})$, as $t\to\ty$. Let $q(t)\in C^3$ be a function comparable of class $1$ to $K t$, $K>0$ in the limit sense, $q''(t)\in o(t^{-3})$, as $t\to\ty$, and $q^{(3)}(t)\in o(t^{-2})$, as $t\to\ty$.

\begin{itemize}

\item[\rm{(i)}] Phase portrait $\Gamma_{xy}=\{(x(t),\dot x(t))\in\eR^2: t\in[t_0,\infty)\}$ of any solution is a spiral near the origin. Phase dimension of any solution of the equation $(\ref{eqalfap})$ is equal to $\dim_{ph}(x)=\frac 2{1+\a}$, for $\a\in(0,1)$.

\item[\rm{(ii)}]  Trajectory $\Gamma $ of the system $(\ref{alfap})$ has box dimension $\dim_B\Gamma =\frac 2{1+\a}$ for $\a\in(0,1)$ and $\gamma\geq\alpha$.

\item[\rm{(iii)}]  Trajectory $\Gamma $ of the system $(\ref{alfap})$ has box dimension $\dim_B\Gamma =2-\frac{\a +\gamma}{1+\gamma}$ for $\a\in(0,1)$ and $0<\gamma<\alpha$.

\item[\rm{(iv)}] Trajectory  $\Gamma $ of the system $(\ref{alfap})$ for $\a >1$ is rectifiable and  $\dim_B\Gamma =1$.
\end{itemize}

\end{theorem}

The graphs of trajectories (\ref{spiral}) for different values of parameter $\alpha$ can be seen in Figures 1--3.

\begin{center}
\includegraphics[width=5cm]{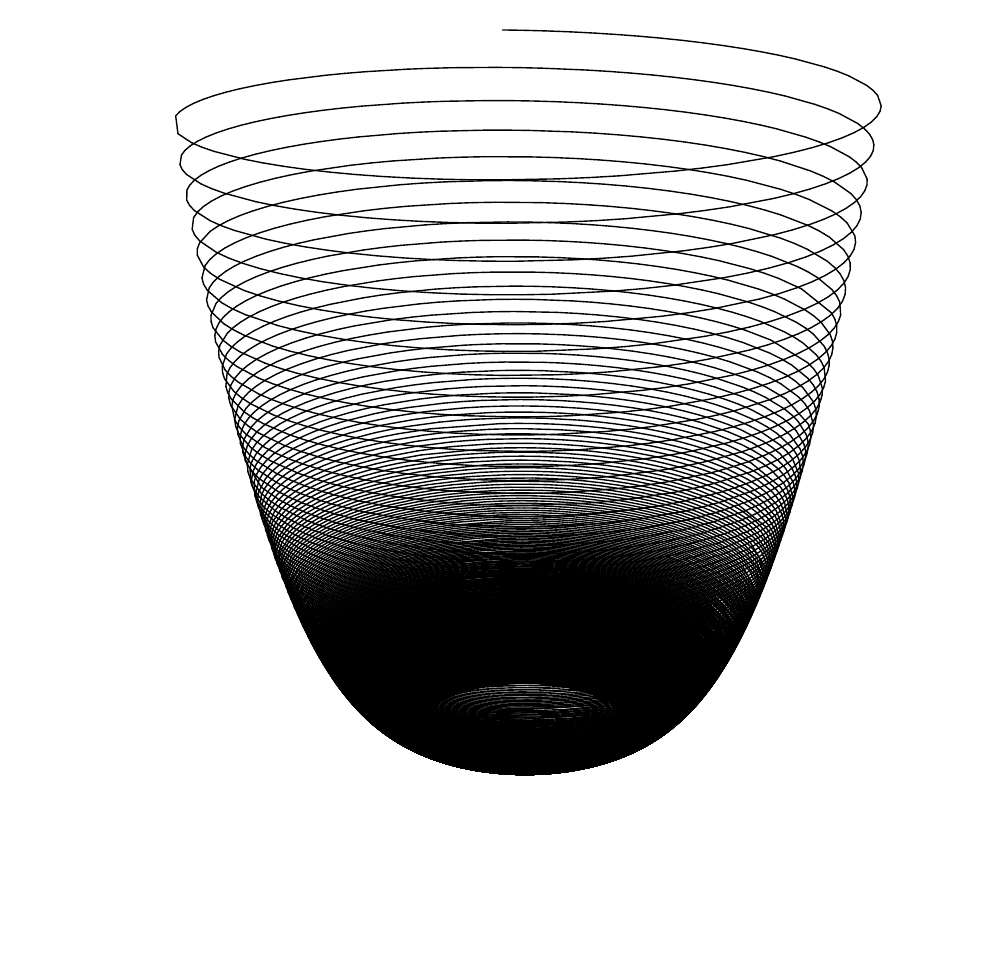} \\
\textbf{Figure 1} System (\ref{alfap}) for $p(t)=t^{-\frac14}$, $q(t)=t$ and $\gamma=1$, Lipschitz case.
\end{center}

\begin{center}
\includegraphics[width=5cm]{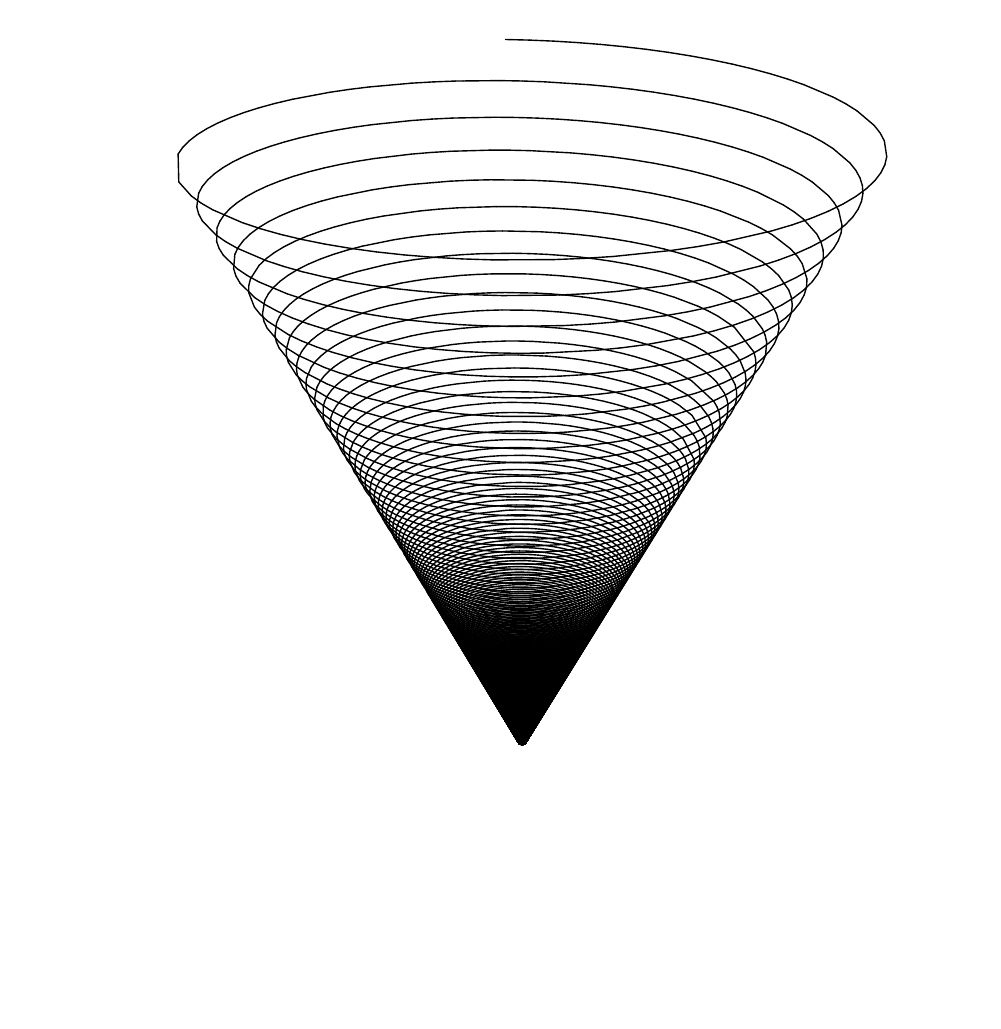} \\
\textbf{Figure 2} System (\ref{alfap}) for $p(t)=t^{-1}$, $q(t)=t$ and $\gamma=1$, Lipschitz case.
\end{center}

\begin{center}
\includegraphics[width=5cm]{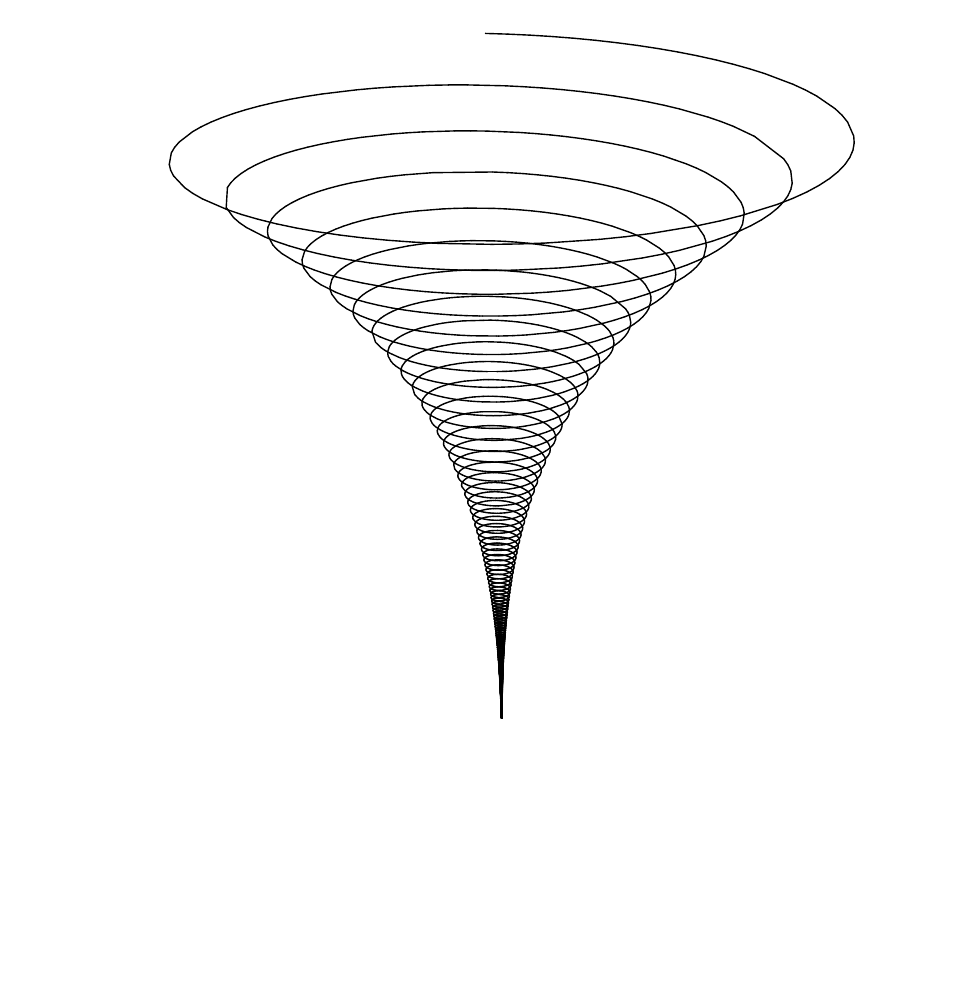} \\
\textbf{Figure 3} System (\ref{alfap}) for $p(t)=t^{-3}$, $q(t)=t$ and $\gamma=1$, H\" older case.
\end{center}

\begin{remark}
In Theorem \ref{phase} the box dimension of the spiral  $\Gamma$  has been computed,  and all values satisfy $\frac{2}{1+\alpha} \leq \dim_B\Gamma< 2-\a$ for $\a\in (0,1)$.
\end{remark}

\begin{remark}
Regarding rectifiability of trajectory  $\Gamma $ of system $(\ref{alfap})$ from Theorem \ref{phase}, the assumptions on functions $p$ and $q$ could be weakened. For instance, for
\bgeqn
p_1(t)&=&t^{-\alpha}\log^k(t) \nonumber\\
p_2(t)&=&t^{-\alpha}\log(\log(\ldots\log(t))),\qquad k\ \mathrm{times} \nonumber\\
q_1(t)&=&t\log^l(t) \nonumber\\
q_2(t)&=&t\log(\log(\ldots\log(t))),\qquad l\ \mathrm{times} \nonumber
\endeqn
where $\alpha>1$ and $k,l\in\mathbb{N}$ if we take $x(t)=p_i(t)\sin(q_j(t))$, $i,j=1,2$, it is easy to see that curve $\Gamma$ is also rectifiable.
If $\alpha\le1$ we expect nonrectifiability and the same box dimension as in the case with no logarithmic terms. This comes from  Remark 9 \cite{zuzu}, saying that  spirals $r=\f^{-\a}(\log\f)^{\b}$, $\f\ge\f_1$, where $\b\ne0$ and $\a\in(0,1)$ have box dimension equal to $d:=2/(1+\a)$ (the same as
for the spiral $r=\f^{-\a}$), but their $d$-dimensional Minkowski content is degenerate. See that degeneracy at Figures 4--6.

We did not prove that all our statements are valid for $p_i, q_i$, $i=1,2$ with logarithmic terms, because in order to do it, we would have to extend theorems from \cite{zuzu} for that cases,  making this article too long. On the other hand, from the dynamical point of view, spirals $r=\f^{-\a}(\log\f)^{\b}$ are not trajectories of vector fields. The Poincar\'e maps or first return maps of foci, which are not weak, have logarithmic terms in the asymptotic expansion. Asymptotics is different in the characteristic directions. These directions could be seen after blowing up when polycycle appears from the focus. The directions with different asymptotic pass through singularities of the polycycle. The logarithmic terms are produced by singularities of the polycycle.
Spiral  $r=\f^{-\a}(\log\f)^{\b}$ has the same asymptotics in all directions, which is a different situation.

\end{remark}

\begin{center}
\includegraphics[width=5cm]{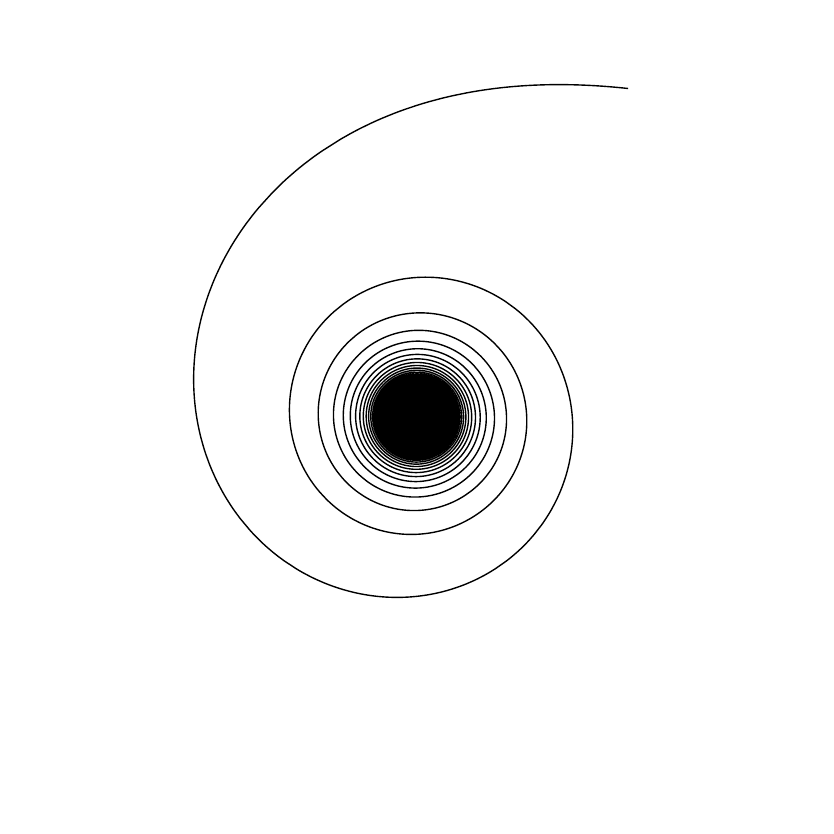} \\
\textbf{Figure 4} Spiral $r=\varphi^{-1/2}$, in polar coordinates.
\end{center}

\begin{center}
\includegraphics[width=5cm]{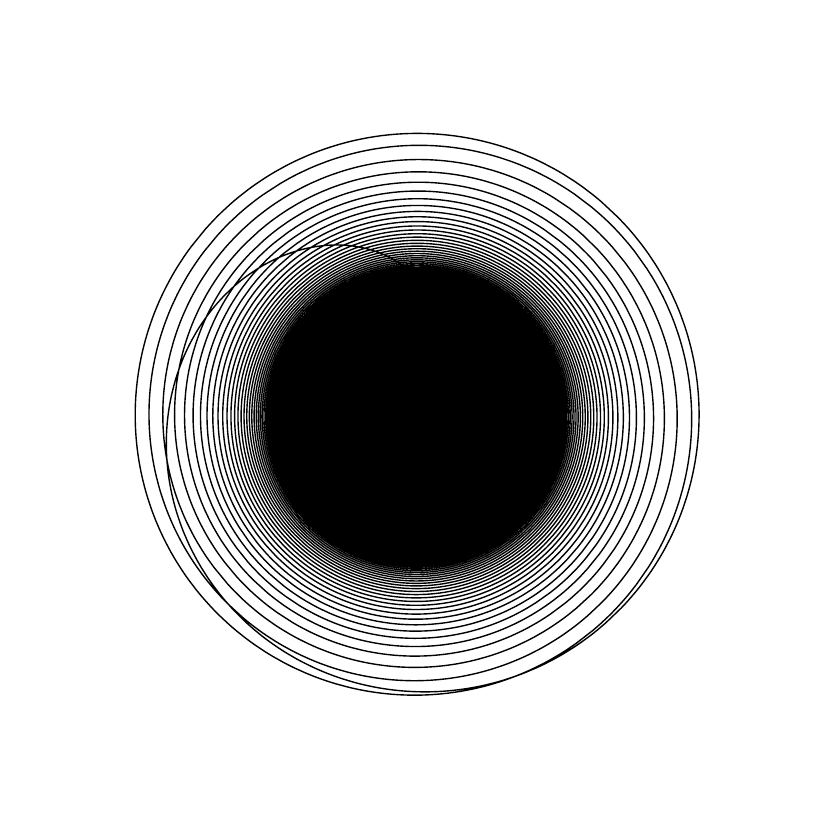} \\
\textbf{Figure 5} Spiral $r=\varphi^{-1/2}\log\varphi$, in polar coordinates.
\end{center}

\begin{center}
\includegraphics[width=5cm]{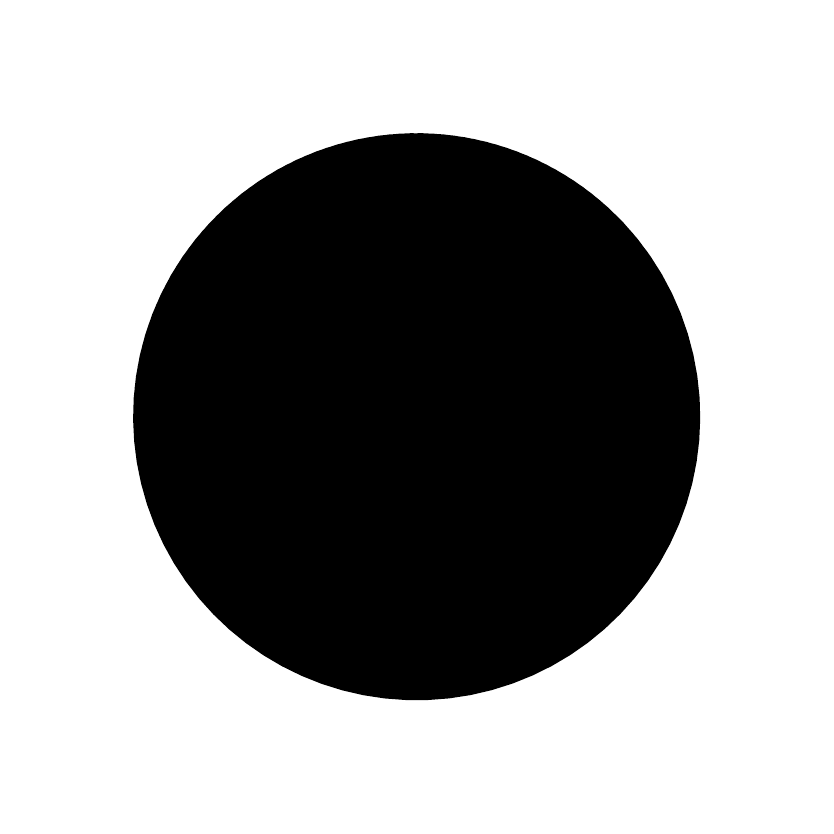} \\
\textbf{Figure 6} Spiral $r=\varphi^{-1/2}\log^2 \varphi$, in polar coordinates.
\end{center}

\begin{remark}
In Introduction and motivation, was briefly explained why we did not take $q(t)\sim t^{\beta}$ for $\beta\neq 1$. Here we would like to show figures concerning that cases.
If $X(\tau)=\tau^{\a}\sin1/\tau^\beta$,  for $\alpha+1\le\beta$ using the described procedure,
we have a planar curve which does not accumulate near the  origin, see Figure 7.

However, if $\alpha+1>\beta$ and $\beta\ne1$, we have spiral converging to zero in ``oscillating'' way, see Figure 8.

Figure 9 shows focus with different asymptotic in the direction of $x$-axes.
\end{remark}

\begin{center}
\includegraphics[width=10cm]{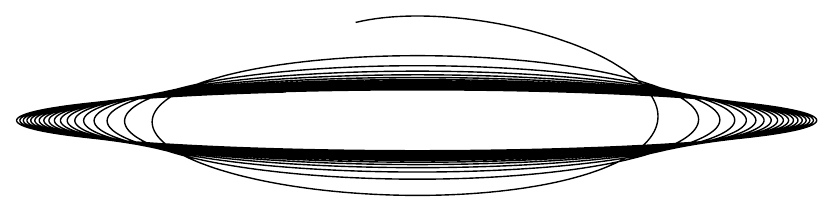} \\
\textbf{Figure 7} Part of unbounded curve $\Gamma_1=\{(x(t),\dot x(t)): t\in[t_0,\infty)\}$, for $x(1/\tau)=X(\tau)=\tau^{1/2}\sin(1/\tau)^{7/4}$, rotated by $\pi/2$ clockwise.
\end{center}

\begin{center}
\includegraphics[width=10cm]{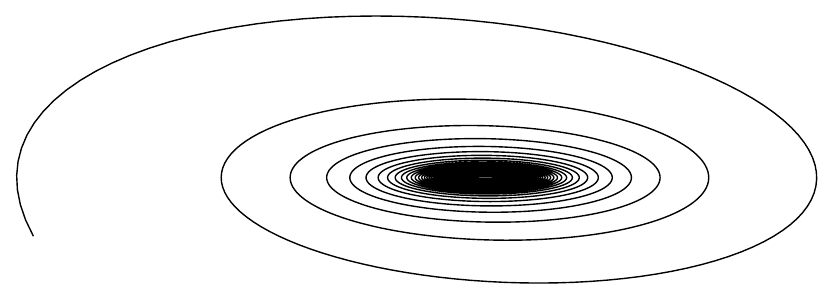} \\
\textbf{Figure 8} Spiral $\Gamma_2=\{(x(t),\dot x(t)): t\in[t_0,\infty)\}$, for $x(1/\tau)=X(\tau)=\tau^{1/2}\sin(1/\tau)^{3/4}$.
\end{center}

\begin{center}
\includegraphics[width=10cm]{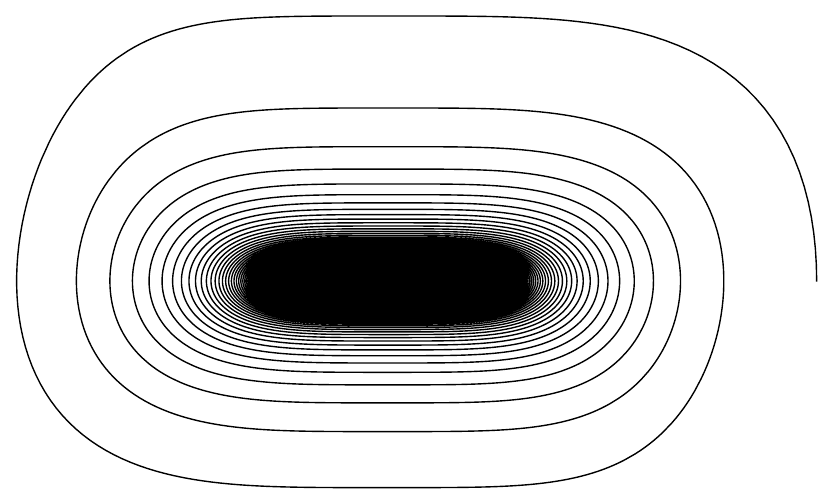} \\
\textbf{Figure 9} Nilpotent focus with characteristic direction along $x$-axes.
\end{center}

\begin{remark}
The system (\ref{alfap}) coincides with the Bessel system of order $\nu$ for $p(t)=t^{-\alpha}$, $\alpha=\nu=1/2$, and $q(t)=t$. The Bessel equation of order $\nu$ has phase dimension equal to $4/3$, for the proof see \cite[Corollary 1]{bessel}. This is a consequence of a fact that the Bessel functions in some sense "behave" like chirps $x(t)=t^{-1/2}\sin\left(t+\theta_0\right)$, $\theta_0\in\eR$, as $t\to\infty$. Although, this background connection is pretty intuitive, the proof is long, complex and technically exhausting.
\end{remark}

The following theorem gives sufficient conditions for rectifiability of a spiral lying into the  H\" olderian  surface $z=g(r)$, $g(r)\simeq r^\b$, $\b>0$. Spiral is called H\"older-focus spiral if it lies in the H\" olderian  surface, and tend to the origin.

\begin{theorem}{\label{rectif}} {\rm (Rectifiability in $\eR^3$)}
 Let $f\:[\f_1,\ty)\to(0,\ty)$,
$\f_1>0$, $f(\f)\simeq \f^{-\a}$,   $|f'(\varphi)|\le C \f^{-\a -1}$,  $\a>1$,  $r=f(\f) $  define a rectifiable spiral. Assume that $g\:(0,f(\f_1))\to(0,\ty)$
is a function of class $C^1$ such that
\bgeq
g(r)\simeq r^\b,\q |g'(r)|\le D r^{\beta -1}, \q \b>0. \nonumber
\endeq
 Let $\C$ be a H\"older-focus spiral defined by $r=f(\f)$, $\f\in[\f_1,\ty)$, $z=g(r)$, then $\C$ is rectifiable spiral.
\end{theorem}
{\it Proof.}
The corresponding parametrization of  spiral $\C$ in Cartesian space coordinates is
\bgeqn\label{system1}\nonumber
\begin{array}{lll}
x&=&f(\f)\cos \f,\\
y&=&f(\f)\sin \f\\
z&=&g(f(\f)).
\end{array}
\endeqn

For the length $l(\Gamma$) of this spiral we have
$$  l(\Gamma)=\int_{\f_1}^{\infty}\sqrt{\dot x^2(\f)+\dot y^2(\f)+\dot z^2(\f)}d\f=$$
$$ \int_{\f_1}^{\infty}\sqrt{f^2(\f)+f'^2(\f)+g'^2(f(\f))f'^2(\f)}d\f\le $$
$$C\int_{\f_1}^{\infty} \sqrt{\f^{-2\a}+\f^{-2\a-2}+\f^{-\a(2\b-2)}\f^{-2\a-2}}d\f=$$
$$C\int_{\f_1}^{\infty} \sqrt{\f^{-2\a}+\f^{-2\a-2}+\f^{-2\a\b-2}}d\f. $$
If $-2\a\b-2\le -2\a$ then
$$ l(\Gamma)\le C\int_{\f_1}^{\infty}\f^{-\a}d\f\ < \ty,$$
and if $-2\a\b-2> -2\a$ then
$$ l(\Gamma)\le C\int_{\f_1}^{\infty}\f^{-\a\b-1}d\f < \ty.$$
\qed

{\it Proof of Theorem \ref{phase}.}
\begin{itemize}
\item[\rm{(i)}] Without the loss of generality we take solution  $x(t)=p(t)\sin q(t)$ of the equation (\ref{eqalfap}). Spiral trajectory $\Gamma$ of system (\ref{alfap}) is defined by (\ref{spiral}). Then $\Gamma_{xy}$ is the projection of $\Gamma$  in $(x,y)-$plane. Using Theorem \ref{gen_wavy}, see the Appendix, we obtain that $\Gamma_{xy}$ is a spiral near the origin and $\dim_B\Gamma_{xy} = \dim_{ph}(x)=\frac 2{1+\a}$.

\item[\rm{(ii)}]
The map  $B\:(x(t),y(t),0)\to(x(t),y(t),z(t))$ is a bi-Lipschitz map near the origin for $\gamma\geq\alpha$, see Lemma \ref{Domagoj}.
It is clear that $\Gamma=B(\Gamma_{xy})$ and it is easy to see that subset $S\subseteq\Gamma$, for which $B$ is not a bi-Lipschitz map, is rectifiable and therefore $\dim_{B}S=1$. The box dimension of set $\Gamma$ is preserved under bi-Lipschitzian mappings and under removing $S\subseteq\Gamma$ such that $\dim_{B}S=1$, see \cite[p.\ 44]{falc}, so it follows form \rm{(i)} that $\dim_B\Gamma=\frac{2}{1+\a}$.

\item[\rm{(iii)}]
Without the loss of generality we take $K=1$. The rest of the proof is similar as in \rm{(ii)}, but using \cite[Theorem 9]{bilip} instead of Lemma \ref{Domagoj}.

\item[\rm{(iv)}]
Without the loss of generality we take $K=1$, because rectifiability is also unaffected by rescaling of spatial variables. Let $r=f(\varphi)$ define curve $\Gamma_{xy}$ in polar coordinates. Notice that $f(\varphi)\simeq\varphi^{-\alpha}$ and $|f'(\varphi)|\leq C\varphi^{-\alpha-1}$, see the proof of Theorem \ref{gen_wavy} from the Appendix. Respecting $r(t)=\sqrt{x(t)^2+\dot{x}(t)^2}$, we take $g(r)$ such that $g(r)\simeq z(t)$ and $g'(r)\in O(z'(t))$, using $z(t)$ from (\ref{spiral}). As $r(t)\simeq t^{-\alpha}$ and $|r'(t)|\leq D_1 t^{-\alpha-1}$, we get $g(r)\simeq r^{\gamma/\alpha}$ and $|g'(r)|\leq D r^{\gamma/\alpha-1}$, so using Theorem \ref{rectif} and the fact that rectifiability is invariant to bi-Lipschitz mapping, as $g(r)\simeq z(t)$, we prove the claim.
\end{itemize}
\qed

\begin{remark}
Notice that in the proof of Theorem \ref{phase} \rm{(iii)}, H\" older case,  we used \cite[Theorem 9]{bilip}, but in the proof of Theorem \ref{phase} \rm{(ii)}, Lipschitz case, we could not use analogous \cite[Theorem 7]{bilip} and we had to devise Lemma \ref{Domagoj}. The reason is behind the fact that assumptions in \cite[Theorem 9]{bilip} about the spiral $r=f(\varphi)$, in polar coordinates, regarding function $f$ being decreasing and $|f'(\varphi)|\simeq\varphi^{-\alpha-1}$, as $t\to\infty$, can be replaced by weaker assumptions. By carefully examining the proof, we see that function $f$ does not have to be decreasing and we can take $|f'(\varphi)|\in O(\varphi^{-\alpha-1})$, as $t\to\infty$. Regardlessly, these assumptions are necessary in \cite[Theorem 7]{bilip}.
\end{remark}

It is interesting to study the Poincar\'e or the first return map associated to a spiral trajectory. The following result is about asymptotics of the Poincar\'e map near focus of the planar spiral from Theorem~\ref{phase}~(i).

\begin{prop}{\rm (Poincar\'e map)} Assume $\Gamma$ is the planar spiral from Theorem~\ref{phase}~(i). Let $P:(0,\varepsilon)\cap\Gamma\to(0,\varepsilon)\cap\Gamma$ be the Poincar\'e map with respect to any axis that passes through the origin.

Then map $P$ has the form $P(r)=r+d(r)$, where $-d(r)\simeq r^{\frac{1}{\a}+1}$ as $r\to 0$.
\end{prop}
{\it Proof.}
Let $\Gamma$ be defined by $r=f(\f)$. Analogously as in the proof of Theorem \ref{gen_wavy}, see \cite[Theorem~4]{cswavy}, it is easy to see that $-d(r)=f(\f)-f(\f+2\pi)\simeq \f^{-\a-1}$  as $\f\to\infty$ and $r\simeq\f^{-\a}$ as $\f\to\infty$. From this follows $-d(r)\simeq r^{\frac{1}{\a}+1}$ as $r\to 0$.
\qed

\medskip

The projection of a solution of system (\ref{alfap}) is a spiral in $(x,y)$-plane. For other two coordinate planes we have the following theorem.

\begin{theorem}{\label{projection}} {\rm (Projections)} Let $p(t)\in C^2$ be a function
 comparable of class $1$ to  power $t^{-\alpha}$, $\alpha>0$, and $p''(t)\in O(t^{-\alpha})$, as $t\to\ty$. Let $q(t)\in C^2$ be a function comparable of class $1$ to $K t$, $K>0$, and $q''(t)\in O(t^{-1})$, as $t\to\ty$.

If $\alpha\in(0,1)$ then projections $G_{xz}$ and $G_{yz}$ of a trajectory (\ref{spiral}), $\gamma>0$ of the system $(\ref{alfap})$ to $(x,z)-$plane and $(y,z)-$plane, respectively, are $(\alpha/\gamma,1/\gamma)$-chirp-like functions, and  $\dim_B G_{xz}=\dim_B G_{yz}=2-\frac{\alpha+\gamma}{1+\gamma}$.
\end{theorem}
{\it Proof.}
Without the loss of generality is $K=1$. Projection $G_{yz}$ is
\bgeqn
Y(z)=y\left(z^{-\frac{1}{\gamma}}\right)=p'\left(z^{-\frac{1}{\gamma}}\right)\sin q\left(z^{-\frac{1}{\gamma}}\right) + p\left(z^{-\frac{1}{\gamma}}\right)q'\left(z^{-\frac{1}{\gamma}}\right)\cos q\left(z^{-\frac{1}{\gamma}}\right) = \nonumber\\
=\sqrt{p'^2\left(z^{-\frac{1}{\gamma}}\right)+p^2\left(z^{-\frac{1}{\gamma}}\right)q'^2\left(z^{-\frac{1}{\gamma}}\right)} \sin\left(z^{-\frac{1}{\gamma}} + \arctan\frac{p\left(z^{-\frac{1}{\gamma}}\right) q'\left(z^{-\frac{1}{\gamma}}\right)}{p'\left(z^{-\frac{1}{\gamma}}\right)}\right). \nonumber
\endeqn
For functions  $P(z)=\sqrt{p'^2\left(z^{-\frac{1}{\gamma}}\right)+ p^2\left(z^{-\frac{1}{\gamma}}\right)q'^2\left(z^{-\frac{1}{\gamma}}\right)}$ and $Q(z)=z^{-\frac{1}{\gamma}} + \arctan\frac{p\left(z^{-\frac{1}{\gamma}}\right) q'\left(z^{-\frac{1}{\gamma}}\right)}{p'\left(z^{-\frac{1}{\gamma}}\right)}$ we have $P(z)\simeq z^{\frac{\alpha}{\gamma}}$,  $P'(z)\simeq z^{\frac{\alpha}{\gamma}-1}$, $Q(z)\simeq z^{-\frac{1}{\gamma}}$, $Q'(z)\simeq z^{-\frac{1}{\gamma}-1}$ as $z\to 0$. So $Y(z)$ is $(\alpha/\gamma,1/\gamma)$-chirp-like function.
To calculate the box dimension of  $G_{yz}$  we apply  \cite[Theorem~5]{cswavy}.

The proof for projection $G_{xz}$ is analogous.
\qed

\begin{remark}
In other words an oscillatory dimension of the solution of (\ref{eqalfap}), under  assumptions of previous theorem concerning $p$ and $q$, is equal to $\dim_{osc}x=\frac{3-\alpha}2$, if  $\a\in (0,1)$.
\end{remark}

\section{Limit cycles}

Limit cycles are interesting object appearing in differential equations.  In particular, we consider a system having
its linear part in Cartesian coordinates
with a conjugate pair $\pm\o i$ of pure imaginary eigenvalues with $\o>0$, and
the third eigenvalue is equal to zero. The corresponding normal form in cylindrical coordinates is:
\bgeqn\label{space}
\dot r&=&a_1 rz+a_2 r^3+a_3rz^2+O(|r,z|)^4\nonumber\\
\dot \f&=&\o+O(|r,z|)^2\\
\dot z&=&b_1r^2+b_2z^2+b_3r^2z+b_4z^3+O(|r,z|)^4,\nonumber
\endeqn
where $a_i$ and $b_i\in\eR$ are coefficients of the system.
Such systems and their bifurcations are treated in Guckenheimer-Holmes \cite[Section 7.4]{gh}. The fold-Hopf bifurcation and  cusp-Hopf bifurcation have been studied in Harlim and Langford \cite{harlim} and the references therein, showing that system (\ref{space}) can exhibit much richer dynamics then singular points and periodic solutions. Notice that in system (\ref{alfap}) there are no limit cycles for any acceptable function $p(t)$. We hypothesize that the limit cycle could be induced by introducing perturbation in the last equation, $\dot{z}=-z^2$.

Here we make a note about box dimension of a spiral trajectory of the simplified system (\ref{space}) at the moment of the birth of limit cycles in $(x,y)$-coordinate plane. In \cite{zuzu}, \cite{zuzulien} we studied planar system consisting of first two equations from (\ref{system}), and  made  fractal analysis of the Hopf bifurcation of the system. We proved that box dimension of a spiral trajectory  becomes nontrivial at the moment of bifurcation. The Hopf bifurcation occurs with box dimension equal to $4/3$, furthermore degenerate Hopf bifurcation or Hopf-Takens bifurcation occurs with the box dimension greater than $4/3$. The more limit cycles have been related to larger box dimension. Analogous results have been showed for discrete systems in \cite{laho}, and applied to continuous systems via Poincar\' e map.
On the other hand in \cite{cras} and \cite{bilip} $3$-dimensional spirals have been studied. Here we consider reduced system
\bgeqn\label{system}
\begin{array}{lll}
\dot r&=&r(r^{2l}+\sum_{i=0}^{l-1}a_ir^{2i})\\
\dot\f&=&1\\
\dot z&=&b_2z^2+\dots+b_nz^n.
\end{array}
\endeqn
First two equations are standard normal form of codimension $l$, where  the Hopf-Takens bifurcation occurs, see \cite{takens}.  The third equation gives us the case where spiral trajectories lie on Lipschitzian or H\" olderian surface, depending on the first exponent. The H\" olderian surface has infinite derivative in the origin, geometrically it is a cusp.

We are interested in the change of the box dimension with respect to the third equation at the moment of birth of limit cycles.
We proved for the standard planar model that the Hopf bifurcation occurs with box dimension equal to $4/3$ and the Hopf-Takens occurs with larger dimensions. Here we prove that on the H\" olderian surface a limit cycle occurs with the box dimension greater than $4/3$.

\begin{theorem}{\label{hopf}} {\rm (Limit cycle)}
Let $l=1$ in the system (\ref{system}) and $b_p<0$ be the first nonzero coefficient in the third equation and $a_0=0$ then  a trajectory $\Gamma$ near the origin has:
\begin{itemize}
\item[\rm{(i)}]
if $2\le p\le 3$ then
$$
\dim_B\Gamma=\frac43,
$$
\item[\rm{(ii)}]
if $p\ge 4$ then
\bgeq\label{holddim}
\dim_B\Gamma=\frac32-\frac{1}{2p}.
\endeq
\end{itemize}
\end{theorem}
{\it Proof.}
Using  \cite[Theorem 9]{zuzu} we get the solution of the first two equations of (\ref{system}), $r\simeq \varphi ^{-1/2}$ having $\dim_B\Gamma_{xy}=\frac43$, where $\Gamma_{xy}$ is orthogonal projection of space trajectory $\Gamma$ to $(x,y)$ plane.
From the third equation we get $z\simeq r^{\frac2{p-1}}$, so for $2\le p\le 3$ we get the Lipschitzian surface, while for $p\ge 4$ surface is H\" olderian. Applying
 \cite[Theorem 7 (a)]{bilip} we obtain $\dim_B\Gamma=4/3$ for $2\le p\le 3$, because the box dimension is invariant for the Lipschitzian case. For the H\" olderian case we apply  \cite[Theorem 9 (a)]{bilip}, where $\alpha=1/2$ and $\beta=2/(p-1)$. So, we get $\dim_B\Gamma=\frac32-\frac{1}{2p}.$

\qed

\begin{remark} The box dimension of a trajectory at the moment of planar Hopf bifurcation is equal to $4/3$, also for $3$-dimensional case with spiral trajectory lying in the  Lipschitzian surface. Situation is different for spiral trajectory contained in the H\" olderian surface, the box dimension of a space spiral trajectory tends to $3/2$. Only one limit cycle could be produced, but dimension increases caused by the H\" olderian behavior near the origin.
Notice that if we apply  formula (\ref{holddim}) obtained for the H\" olderian case, to  the Lipschitzian case $p=3$ we will get correct result $4/3$.
For $l>1$ degenerate Hopf bifurcation or  Hopf-Takens  bifurcation appears, where $l$ limit cycles could be born, and  the box dimension of the space spiral trajectory is equal to $\dim_B\Gamma=\frac{(4l-1)p-2l+1}{2lp}$ using the same arguments.
\end{remark}

\appendix
\section{Auxiliary results}

\begin{theorem}{\label{gen_wavy}} {\rm (Generalization of \cite[Theorem~4]{cswavy})}
Let $p(t)\in C^3$ be a function comparable of class $2$ to  power $t^{-\alpha}$, $\alpha>0$, in the limit sense, and $p^{(3)}(t)\in O(t^{-\alpha-3})$, as $t\to\ty$. Let $q(t)\in C^3$ be a function comparable of class $1$ to $K t$, $K>0$ in the limit sense, $q''(t)\in o(t^{-3})$, as $t\to\ty$, and $q^{(3)}(t)\in o(t^{-2})$, as $t\to\ty$.

Define $x(t)=p(t) \sin q(t)$ and continuous function $\varphi (t)$ by $\tan \f(t)=\frac{\dot{x}(t)}{x(t)}$.

\begin{itemize}
\item[\rm{(i)}] If $\alpha\in(0,1)$ then the planar curve $\Gamma:=\{(x(t),\dot x(t))\in\eR^2: t\in[t_0,\infty)\}$ is a spiral $r=f(\f)$, $\f\in(-\infty,-\phi_0]$, near the origin, and
    $$
    \dim_{ph}(x):=\dim_B\Gamma=\frac{2}{1+\alpha} .
    $$
\item[\rm{(ii)}] If $\alpha>1$ then the planar curve $\Gamma$ is a rectifiable spiral near the origin.
\end{itemize}

\end{theorem}
{\it Proof.}
After substitution of time variable by $u=t/K$ and respective rescaling of the $x$ and $y$ axes, we continue assuming $K=1$. The rest of the proof is analogous to the proof of \cite[Theorem~4]{cswavy}, but carefully taking care about the more general conditions on $q$. Rescaling of the $x$ and $y$ axes in the plane is a bi-Lipschitz map, so the box dimension remains preserved.
\qed

\bigskip

\bibliographystyle{plain}
\bibliography{bibliografija}

\end{document}